\documentclass[3p]{elsarticle}

\usepackage{ulem}
\usepackage{color} 
\usepackage{psfrag} 
\usepackage{subfig}
\usepackage{hyperref}
\usepackage{csvsimple}
\usepackage{amsmath}
\usepackage{amsthm}
\usepackage{graphics}
 \usepackage{pstool} 
\usepackage
{todonotes}

\usepackage{pgfplots}
\usepackage{pgfplotstable}
\usepackage{filecontents}
\usepackage{pgf}
\usepgfplotslibrary{colorbrewer}

\usepackage{tikz}
\usetikzlibrary{decorations.pathmorphing}
\usepackage{tikzscale}
\usepackage{tikz-3dplot}
\usepackage{graphicx}

\usetikzlibrary{3d}
\usetikzlibrary{spy}
\usetikzlibrary{snakes,arrows,shapes}    
\usetikzlibrary{spy}
\usetikzlibrary{backgrounds}  
\usetikzlibrary{decorations}
\usetikzlibrary{positioning}
\usetikzlibrary{patterns}
\usetikzlibrary{calc} 
\usetikzlibrary{external} 

\usepackage{booktabs}   
\usepackage{makecell} 
\usepackage{colortbl}   
\usepackage{xspace}
\usepackage{color}     
\usepackage{setspace}   
 
\usepackage{soul}
\usepackage{ulem}
\usepackage{currfile}
\usepackage{import}

\usepackage{floatrow}
\usepackage{multicol, blindtext}
\newfloatcommand{capbtabbox}{table}[][\FBwidth]

\pgfplotsset{compat=1.7}

\newcommand{\findmax}[3]{
    \pgfplotstablesort[sort key={#2},sort cmp={float >}]{\sorted}{#1}%
    \pgfplotstablegetelem{0}{#2}\of{\sorted}%
    \let #3=\pgfplotsretval%
}

\pgfplotscreateplotcyclelist{myCycleListBlackAndWhite}
{%
  {black, mark=o},
  {black, mark=square},
  {black, mark=triangle},
  {black, mark=diamond}, 
  {black, mark=pentagon}, 
  {black, mark=*},
  {black, mark=square*},
  {black, mark=triangle*},
  {black, mark=diamond*}, 
  {black, mark=pentagon*}, 
}

\definecolor{darkgreen}{rgb}{0,0.4,0} 
\definecolor{darkbrown}{rgb}{0.5, 0.396, 0.09}

\definecolor{c1}{rgb}{0.0, 0.4196078431372549, 0.6431372549019608}
\definecolor{c2}{rgb}{1.0, 0.5019607843137255, 0.054901960784313725}
\definecolor{c3}{rgb}{0.6705882352941176, 0.6705882352941176,
0.6705882352941176} \definecolor{c}{rgb}{0.34901960784313724, 0.34901960784313724, 0.34901960784313724}
\definecolor{c4}{rgb}{0.37254901960784315, 0.6196078431372549,
0.8196078431372549} \definecolor{c}{rgb}{0.7843137254901961, 0.3215686274509804, 0.0}
\definecolor{c5}{rgb}{0.5372549019607843, 0.5372549019607843,
0.5372549019607843} \definecolor{c}{rgb}{0.6352941176470588, 0.7843137254901961, 0.9254901960784314}
\definecolor{c6}{rgb}{1.0, 0.7372549019607844, 0.4745098039215686}
\definecolor{c7}{rgb}{0.8117647058823529, 0.8117647058823529,
0.8117647058823529}

\pgfplotscreateplotcyclelist{myCycleListColor}
{%
  {blue, mark=o},
  {red, mark=square},
  {darkbrown, mark=triangle},
  {darkgreen, mark=diamond},
  {black, mark=pentagon},
  {blue, mark=*},
  {red, mark=square*},
  {darkbrown, mark=triangle*},
  {darkgreen, mark=diamond*},
  {black, mark=pentagon*},
}

\pgfplotscreateplotcyclelist{myCycleListColorNoMarkers}
{
  {black},
  {blue, dashed},
  {red, dotted},
  {darkbrown, loosely dashed},
  {darkgreen, loosely dotted},
}

\pgfplotsset{every axis/.append style= 
              {
                font=\small,
                mark size=2,
                line width = 0.1,
                legend style={font=\small, mark size=3, draw=none, fill=none},
                legend cell align=left,
                cycle list name=myCycleListColor,
              }
            }
            
\pgfplotstableset
{
  every odd row/.style={before row={\rowcolor[gray]{0.9}}},
  every head row/.style=
  {
    before row=
    {%
      \toprule
    },
    after row=\midrule
  },
  every last row/.style=
  {
    after row=\bottomrule
  }
}

\newif\ifdrawboundingbox

\usetikzlibrary{external} 
\tikzset{external/system call={pdflatex \tikzexternalcheckshellescape
-halt-on-error -interaction=batchmode -jobname "\image" "\texsource"}} 

\pgfkeys
{ 
  /pgf/number format/.cd,
  fixed,
  set thousands separator={\,},
  1000 sep in fractionals,
}

\newcolumntype{C}[1]{>{\centering\arraybackslash}m{#1}}
\newcolumntype{R}[1]{>{\raggedright\arraybackslash}m{#1}}
\newcolumntype{L}[1]{>{\raggedleft\arraybackslash}m{#1}}

\newcommand{\delete}[1]{\xspace}

\drawboundingboxtrue
\drawboundingboxfalse 

\usepackage{amsmath,amsfonts}

\usepackage{etoolbox}
\ifdef{\R}{
	\renewcommand{\R}{\mathbb{R}}
}{
\newcommand{\R}{\mathbb{R}}
}
\ifdef{\N}{\renewcommand{\N}{\mathbb{N}}}{\newcommand{\N}{\mathbb{N}}}
\ifdef{\C}{\renewcommand{\C}{\mathbb{C}}}{\newcommand{\C}{\mathbb{C}}}
\ifdef{\Z}{\renewcommand{\Z}{\mathbb{Z}}}{\newcommand{\Z}{\mathbb{Z}}}

\newcommand{\lnorm}[1]{{\left\vert\kern-0.25ex\left\vert\kern-0.25ex\left\vert #1 
    \right\vert\kern-0.25ex\right\vert\kern-0.25ex\right\vert}}

\newcommand{\meinlemma}{Lemma \stepcounter{smeinTheorems}\the\value{smeinTheorems}}
\newcommand{\meinsatz}{Satz \stepcounter{smeinTheorems}\the\value{smeinTheorems}}
\newcommand{\meinedefinition}{Def.\ \stepcounter{smeinTheorems}\the\value{smeinTheorems}}

\newcounter{meineAussage} 
\newcommand{\opnorm}[1]{{\left\vert\kern-0.25ex\left\vert\kern-0.25ex\left\vert #1 
    \right\vert\kern-0.25ex\right\vert\kern-0.25ex\right\vert}}

\usepackage{mathtools}
\usepackage{multicol}
\usepackage{amssymb}
\usepackage{amsbsy}

\usepackage{svg}
\usepackage{algorithm2e}
\usepackage{algpseudocode}
\usepackage{gensymb}
\usepackage{todonotes}
\usepackage[capitalise,noabbrev,sort&compress]{cleveref}
\usepackage{soul,xcolor}
\usepackage{gensymb}
\usepackage[T1]{fontenc}
\usepackage{accents}
\newcommand{\timeInterval}{\lbrack t_{0},\,t_{end}\rbrack}

\newcommand{\theParameters}{$\rho$, $c$, and $\kappa$ }

\begin{document}

\begin{frontmatter}

\title{A spatiotemporal two-level method for high-fidelity thermal analysis of laser powder bed fusion}
\author[GSSIAddress]{Alex Viguerie}

\author[PaviaAddress]{Massimo Carraturo\corref{mycorrespondingauthor}}
  \cortext[mycorrespondingauthor]{Corresponding author}
  \ead{massimo.carraturo@unipv.it}

\author[PaviaAddress]{Alessandro Reali}

\author[PaviaAddress]{Ferdinando Auricchio}

\address[GSSIAddress]{Gran Sasso Science Institute, 
viale F. Crispi 7, 67100 L`Aquila, Italy} 

\address[PaviaAddress]{Department of Civil Engineering and Architecture, 
University of Pavia,
via Ferrata 3, 27100 Pavia, Italy}

\newcommand{\publicationDate}{\today}


\vspace{-1.5cm} 
\hrule 

\begin{abstract}
Numerical simulation of the laser powder bed fusion (LPBF) procedure for additive manufacturing (AM) is difficult due to the presence of multiple scales in both time and space, ranging from the part scale (order of millimeters/seconds) to the powder scale (order of microns/milliseconds). This difficulty is compounded by the fact that the regions of small-scale behavior are not fixed, but change in time as the geometry is produced. While much work in recent years has been focused on resolving the problem of multiple scales in space, there has been less work done on multiscale approaches for the temporal discretization of LPBF problems. In the present contribution, we extend on a previously introduced two-level method in space by combining it with a multiscale time integration method. The unique transfer of information through the transmission conditions allows for interaction between the space and time scales while reducing computational costs. At the same time, all of the advantages of the two-level method in space (namely its geometrical flexibility and the ease in which one may deploy structured, uniform meshes) remain intact. Adopting the proposed multiscale time integration scheme, we observe a computational speed-up by a factor $\times 2.44$ compared to the same two-level approach with uniform time integration, when simulating a laser source traveling on a bare plate of nickel-based superalloy material following an alternating scan path of fifty laser tracks.
\end{abstract}
 

 

\end{frontmatter}

\sloppy

\section{Introduction}\label{sec:introduction}


Additive manufacturing (AM) technologies have grown significantly in recent years due to their ability to produce designs difficult or even impossible to obtain with standard manufacturing techniques \cite{king_laser_2015}. Though multiple technologies for AM exist, laser powder bed fusion (LPBF) is among the most commonly used, and is the focus of the current work. LPBF starts via the deposition of a thin powder layer over a build plate, which is then selectively heated by a laser (in accordance with the desired geometry), melting the powder into liquid. Upon cooling, this liquid becomes solid and a new layer of powder can then be deposited. Such a process is repeated until the desired object has been produced. The rapid melting-solidification cycles occurring during the process induce residual stresses leading to undesirable features and design flaws and eventually even to process failures \cite{herzog2016additive,ghosh_single-track_2018}. Reliable numerical simulation of LPBF processes can help to prevent residual stresses and thus there is a huge interest in both the academic and industrial communities in developing and implementing effective numerical schemes to achieve such a goal.

\par However, the simulation of LPBF is extremely difficult, owing to wide range of relevant spatiotemporal scales in the problem. In fact, the laser melting-solidification cycle occurs at the powder scale (on the order of microns in space and milliseconds in time); at the same time, larger-scale features, on coarser space-time scales (millimeters and seconds, respectively) are nonetheless important and cannot be neglected, even when we limit the investigation to quantities at the melt pool scale, i.e., melt pool morphology, cooling rate, etc. As there is interaction between these scales, and their behaviors influence one another, decoupling solutions is also not an option. Further complicating the problem, regions of large- and small-scale behavior change in time as the laser moves in order to produce the desired geometry. Thus, efficient simulation of LPBF processes requires methods that are capable of resolving both fine- and coarse-scale behavior in both space and time while remaining geometrically flexible. 
 
 \par Several approaches have been employed to resolve the multiple scales of LPBF problems in space. A common method to address this difficulty is through the use of refinement-and-derefinement algorithms \cite{Patil2015, carraturo2019a, kollmannsberger_hierarchical_2018, LI2019100903, baiges2020adaptive}. Such approaches, however, may be difficult to implement, expensive, and lead to numerical conditioning problems. As an alternative, in \cite{viguerie2020fat}, a \textit{two-level} method was proposed, in which two different meshes are used: a \textit{global} mesh for the coarse-scale behavior and a \textit{local} mesh for the fine scale behavior; this approach retains the advantages of refinement-and-derefinement while allowing for the use of uniform meshes. In \cite{viguerie2021numerical}, the two-level method was shown to be a viable choice also for LPBF simulations.  
 
\par However, there has been comparatively less work on separating the multiple scales in time present in LPBF processes. A recent work by \citet{hodge2021towards} introduces a framework for multiscale (or multirate) time integration for additive problems. In this work, the different (spatial) degrees of freedom are partitioned according to whether they represent `fast' or `slow' dynamics, and the two scales interact with each other via a predictor-corrector type scheme. Hodge acknowledges the possibility that the spatial regions showing `fast' and `slow' behavior may change in time and suggests different ways of overcoming such practical difficulties, including heuristic, physics-informed, and algorithmic approaches to identifying the time scale of the different sets of degrees of freedom. In \cite{hodge2021towards}, it  was then demonstrated that multirate time integration can deliver results similar to fine-scale uniform time integration at a small fraction of the computational cost. 


\par In \cite{SOLDNER20192183}, Soldner and Mergheim employ a domain-decomposition approach in which the different time scales were separated by physical or similar arguments, and then solve in a discountinuous-Galerkin type framework in time, in which penalty terms along the interface are included in order to resolve the possible time discontinuities. The interface between the regions of `fast' and `slow' dynamics is then resolved explicitly, leading to a critical time step size. A \textit{merge-and-split} scheme was also introduced to demonstrate how this approach could be easily adapted in time as the fast dynamics regions change. It was shown through a series of numerical examples to be an effective approach in resolving thermal and additive-type problems.


\par \citet{CHENG2021113825} present a sub-structuring domain decomposition method for transient heat transfer problems. Starting from the work of \citet{roux2009domain}, they employ a Dirichlet-Robin iteration method to solve time-dependent problems, showing a significant reduction in computational costs, unconditional stability, and rapid convergence for a 2D model of a single laser track AM problem setup.

\par Outside of the specific application to additive manufacturing, multi-scale time integration was also the subject of some research in the 1970s, with various works by Hughes \cite{hughes1978implicit, miranda1989improved} and Belytschko \cite{belytschko1979mixed, belytschko1978stability, liu1984partitioned} being of particular note. These works also concentrate with spatially-partitioned multirate integration schemes. The different regions are often integrated with different schemes; such works focus on \textit{implicit}-\textit{explicit}, \textit{explicit}-\textit{explicit}, and \textit{implicit}-\textit{implicit} methods. More recent works by Sandu et. al \cite{sandu2019class, gunther2016multirate} focus on Runge-Kutta type methods in a more general setting. The term `multirate' was introduced by Sandu et. al in \cite{gunther2016multirate}, and then used in \cite{hodge2021towards}, and will also be used here, interchangeably with `multiscale.'

\par In the present work, we introduce a method that combines the two-level method for spatial scale separation \cite{viguerie2020fat, viguerie2021numerical} with a multirate time integration scheme for the resolution of the different scales in both space and time for the LPBF problem. The described approach differs from those shown in \cite{hodge2021towards, SOLDNER20192183} in that it employs multiple meshes, and by the mechanism through which information is passed between the local and global scales. In the present work, we use a line (for two dimensions) or a surface (for three dimensions) integral formulation inspired by the fat-boundary method \cite{maury2001fat, bertoluzza2011analysis,bertoluzza2005fat} to exchange information between local and global problems. We note as well that the use of distinct meshes simplifies the problem of identifying the different regions of time integration; we may define one time scale as belonging to one mesh. This approach also grants significant geometric flexibility; as there is no need to explicitly update the degrees of freedom requiring a finer temporal resolution as the laser moves its position in time.
\par The article is outlined as follows. We begin by introducing the relevant physical models. We briefly discuss the two-level method in space and recall its formulation and basics regarding its implementation. We then introduce and describe the two-level time integration scheme, emphasizing its natural connection with the corresponding spatial scheme. We finally validate with a series of numerical experiments, establishing both the accuracy and computational viability of the introduced method, demonstrating its ability to capture realistic problem features across multiple space and time scales.


%



\section{Governing equations}\label{sec:governingEquations}
In this section we present the model equations used herein to describe the laser powder bed fusion process. An effective model must incorporate many different physical phenomena, in particular nonlinear heat transfer, the phase change between powder, liquid, and solid material, and the resulting influence of latent heat due to material phase change. We consider herein  the two-phase latent-heat model of \cite{kollmannsberger2019}. \par  The adopted model consists of two coupled parts: a \textit{thermal model} describing the evolution of the temperature field in the material and a \textit{material model} describing the change of the physical material properties as the temperature changes. We will present the different components separately.
\subsection{Thermal model}
The thermal model is described through an unsteady partial differential equation in a domain $\Omega_+$ over a time interval $\lbrack t_{0},\,t_{end}\rbrack$ in terms of a temperature $T$:
\begin{equation}\label{thermalEqn}\begin{alignedat}{2}
\rho c \frac{\partial T}{\partial t} + \chi \rho \frac{\partial f_{pc}}{\partial t} - \nabla \cdot \left(\kappa \nabla T\right)   &= Q \qquad &&\text{in } \Omega_+ \times \timeInterval \\
T &= T_D \qquad &&\text{on } \Gamma_D \times \timeInterval  \\
T &= T_0 \qquad &&\text{in } \Omega_+, \,\,t=t_0,
\end{alignedat}\end{equation}
where $Q$ is a heat source term while $\rho$, $c$, $\kappa$, and $\chi$ denote the density, specific heat capacity, thermal conductivity, and latent heat respectively. 
 In general, the thermal parameters \theParameters are functions of the temperature $T$; this dependence is understood if not denoted explicitly and will be detailed in the following subsection. The term $f_{pc}$ denotes the \textit{phase change} function, and models the effect of latent heat released by such changes. $f_{pc}$ itself is chosen to have a sigmoid behavior. The specific definition of $f_{pc}$ is problem- and material-dependent and its exact definitions will be specified as necessary.


\subsection{Material Model}
\begin{figure}[h!]
	\centering
	\subfloat[Thermal conductivity.\label{subfig:localBCs}]
	{
		\includegraphics[width=0.5\textwidth]{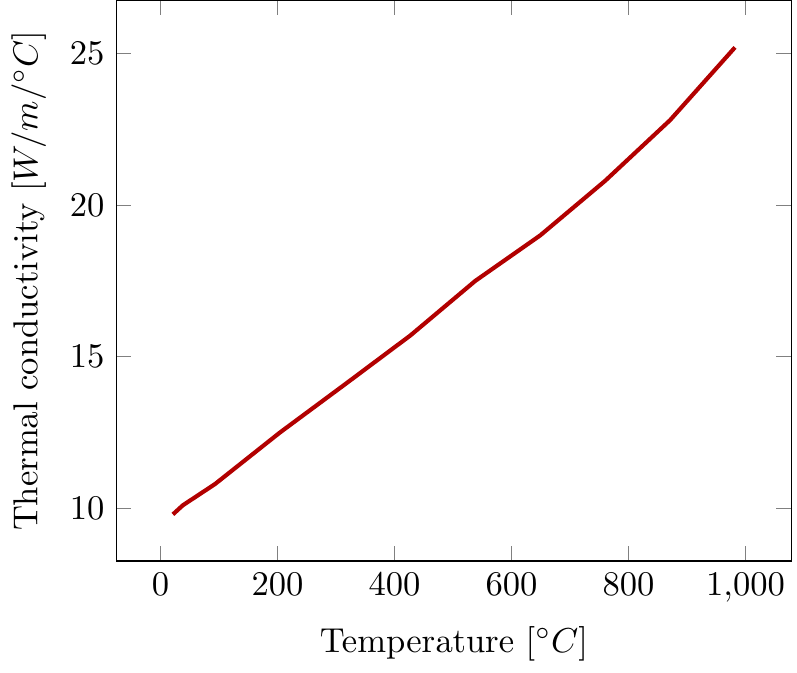}
	}
	\\
	\subfloat[Specific heat capacity.\label{subfig:globalBCs}]
	{
		\includegraphics[width=0.5\textwidth]{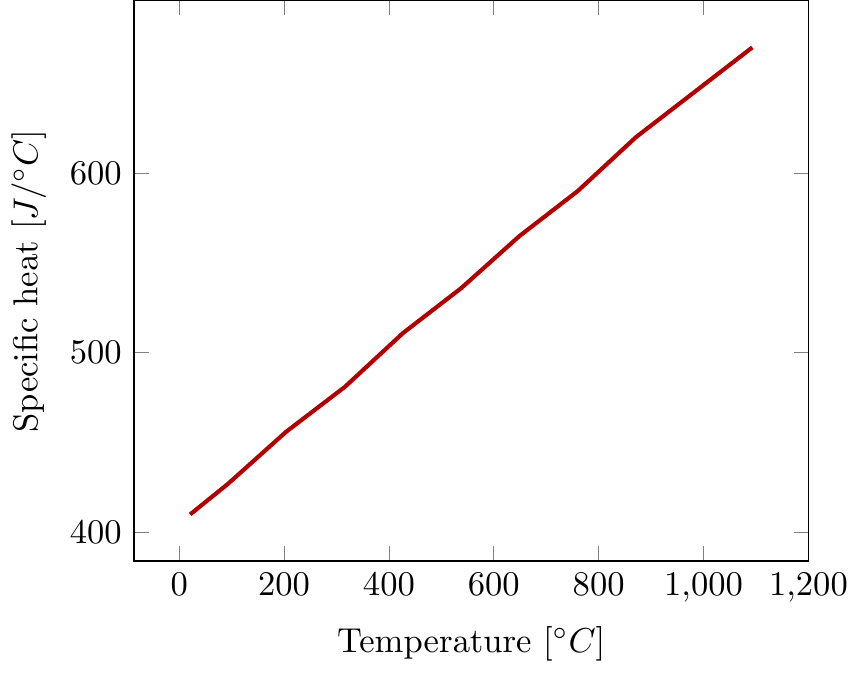}
	}    
  \caption{IN625 Temperature-dependent material properties taken from ~\citep{specialmetals}.\label{fig:IN625TempMat}}
\end{figure}

The nickel-based superalloy IN625 has been considered as reference material in the present contribution due to its widespread adoption in several LPBF applications and its well-known material properties ~\cite{specialmetals}.
Values of the temperature dependent specific heat capacity and thermal conductivity are reported in~\cref{fig:IN625TempMat}, whereas parameter values above the maximum temperature in the plot have been set to constant. In the present work, we assume a constant density value $\rho=8440$ [kg/m$^3$].

\section{The two-level method in space}\label{ssec:TwoLevelMethod}

The two-level method is based on a re-formulation of the heat transfer equation \eqref{thermalEqn} as two coupled problems referred in the following as the \textit{local} and the \textit{global} problem.
Let us consider a global domain $\Omega_+$ and a local domain $\Omega_-\subset \Omega_+$, such that, following the notation introduced in \cite{viguerie2020fat}, $T^+$ and $T^-$ are the local and the global temperature fields respectively, whereas the thermal conductivity $\kappa$ is defined as:
\begin{equation}
	\kappa := 
	\begin{cases}
		\kappa_+ & \text{in } \Omega_+ \setminus \Omega_- \\
		\kappa_- & \text{in } \Omega_-
	\end{cases}
\end{equation}
with $\rho_+$, $\rho_-$, $c_+$, and $c_-$ defined analogously. Note that these may be functions of time, space, and temperature in general, and this dependence is assumed if not explicitly denoted. We denote as $\gamma$ the \textit{interface} between local and global problems across which information is exchanged. \par For $\eta \in H^{-1/2}\left(\gamma \right)$, we define $\eta \xi_{\gamma} \in H^{-1}\left(\Omega_+\right) = H^1(\Omega_+)'$ as the linear functional such that:
  \begin{align}\label{deltaDef}
  	\int_{\Omega_+} \left(\eta \xi_{\gamma}\right)w = \int_{\gamma} \eta w \quad\quad \forall w \in H_{\Gamma_D}^1\left(\Omega_+\right),
  \end{align} 
  \noindent following the notation used in \cite{maury2001fat, bertoluzza2011analysis, bertoluzza2005fat, viguerie2020fat, viguerie2021numerical}.
\par We then define the \textit{global problem} as:
 \begin{align}\begin{split}\label{globalPbm}\MoveEqLeft[2]
c_+\rho_+\frac{\partial T^+}{\partial t} + \left(\frac{c_- \rho_- \kappa_+}{\kappa_-} - c_+\rho_+\right)\frac{\partial T^-}{\partial t}\big|_{\Omega_-} \\ \quad{}&- \nabla \cdot \left(\kappa_+ \nabla T^+\right) = Q\big|_{\Omega_+} + \frac{\kappa_+}{\kappa_-} Q\big|_{\Omega_-} \\ &\quad{}+ \frac{\kappa_+}{\kappa_-} \nabla T^- \cdot \nabla\kappa_-\big|_{\Omega_-} - \nabla T^- \cdot \nabla \kappa_+\big|_{\Omega_-} \\ &+ \left(\kappa_+ - \kappa_-\right)\frac{\partial T^-}{\partial \boldsymbol{n} } \xi_{\gamma}\hspace{1.cm}\text{  in  }\Omega_+ \times \lbrack t_0,\,t_{end}\rbrack \\
T^+ &= T_0 \hspace{3.3cm}\text{  on  } \Gamma_D \\
\kappa_+ \frac{\partial T^+}{\partial \boldsymbol{n}} &=0 \hspace{3.5cm}\text{  on  } \Gamma_N,
 \end{split}\end{align}
 with the corresponding \textit{local problem} given by:
 \begin{equation}\label{localProblem}\begin{alignedat}{2}
 	&c_-\rho_- \frac{\partial T^-}{\partial t} + \chi \rho_- \frac{\partial f_{pc}}{\partial t} - \nabla \cdot \left(\kappa_- \nabla T^- \right) = Q \qquad &&\text{  in  } \Omega_- \times \lbrack t_0,\,t_{end}\rbrack \\
 	&T^- = T^+ \qquad &&\text{  on  } \gamma_D \\
 	&\kappa_- \frac{\partial T^- }{\partial\boldsymbol{n} } = 0 \qquad &&\text{  on  } \Gamma_N .
\end{alignedat}\end{equation}
It was shown in \cite{viguerie2020fat} that the problems (\ref{globalPbm}, \ref{localProblem}) are equivalent to (\ref{thermalEqn}) for non-constant thermal parameters. Note here that we only consider the effects of latent heat in the local problem. In practice, one obtains the two-level method by solving the local and global problems iteratively, in a (possibly relaxed) scheme, until convergence is reached.
\par One may instead employ an alternate formulation of the global problem, which reads:
 \begin{align}\begin{split}\label{globalPbmAlternate}
&c_+\rho_+\frac{\partial T^+}{\partial t}  - \nabla \cdot \left(\kappa_+ \nabla T^+\right) \\ &= Q+ \left(\kappa_+ - \kappa_-\right)\frac{\partial T^-}{\partial \boldsymbol{n} } \xi_{\gamma} \hspace{.5cm}\text{  in  }\Omega_+ \times \lbrack t_0,\,t_{end}\rbrack \\
T^+ &= T_0 \hspace{3.6cm}\text{  on  } \Gamma_D \\
\kappa_+ \frac{\partial T^+}{\partial \boldsymbol{n}} &=0 \hspace{3.75cm}\text{  on  } \Gamma_N.
 \end{split}\end{align}
This may offer significant advantages in terms of computational time and conditioning. As the derivatives of $\kappa$ are no longer required, there is no need to store $\kappa$ in a way that allows for numerical differentiation, and the numerical integration during the assembly phase is greatly simplified. 
\par The disadvantage of \eqref{globalPbmAlternate} is a mild loss of physical consistency in the following sense: if $T^*$ is the solution of the full coupled problem \eqref{thermalEqn}, $T^+ \in X_+(\Omega_+)$ a weak solution of \eqref{globalPbmAlternate} in a suitable function space $X_+$, and $T^- \in X_-(\Omega_-)$ a weak solution of \eqref{localProblem} in a suitable function space $X_-$, then we have that:
\begin{align}
(T^*-T^-,\varphi)_{\Omega_-} &= 0\qquad \forall \varphi \in X_-, \\
(T^* - T^+,\phi)_{\Omega_+ \setminus \Omega_- } &= 0 \qquad \forall \phi \in X_+.
\end{align}
However, the more stringent consistency condition:
\begin{align}\label{superConsistency}
(T^* - T^+,\phi)_{\Omega_+ \cap \Omega_-} &= 0 \qquad \forall \phi \in X_+,
\end{align}
which holds for a solution of \eqref{globalPbm}, does not hold in general for a solution of \eqref{globalPbmAlternate}. In practice, this error is not large, and a more rigorous investigation of this consistency error will be the subject of future work. We note that for many problems, the scale separation between spaces spanned by $X_-$ and $X_+$ is sufficiently large such that relaxing the condition \eqref{superConsistency} does not lead to appreciable errors.
\subsection{Local and global problem boundary conditions}\label{ssec:thermalBC}
\begin{figure}[h!]
	\centering
	\subfloat[Global domain BCs.\label{subfig:globalBCs}]
	{
		\includegraphics[width=0.75\textwidth]{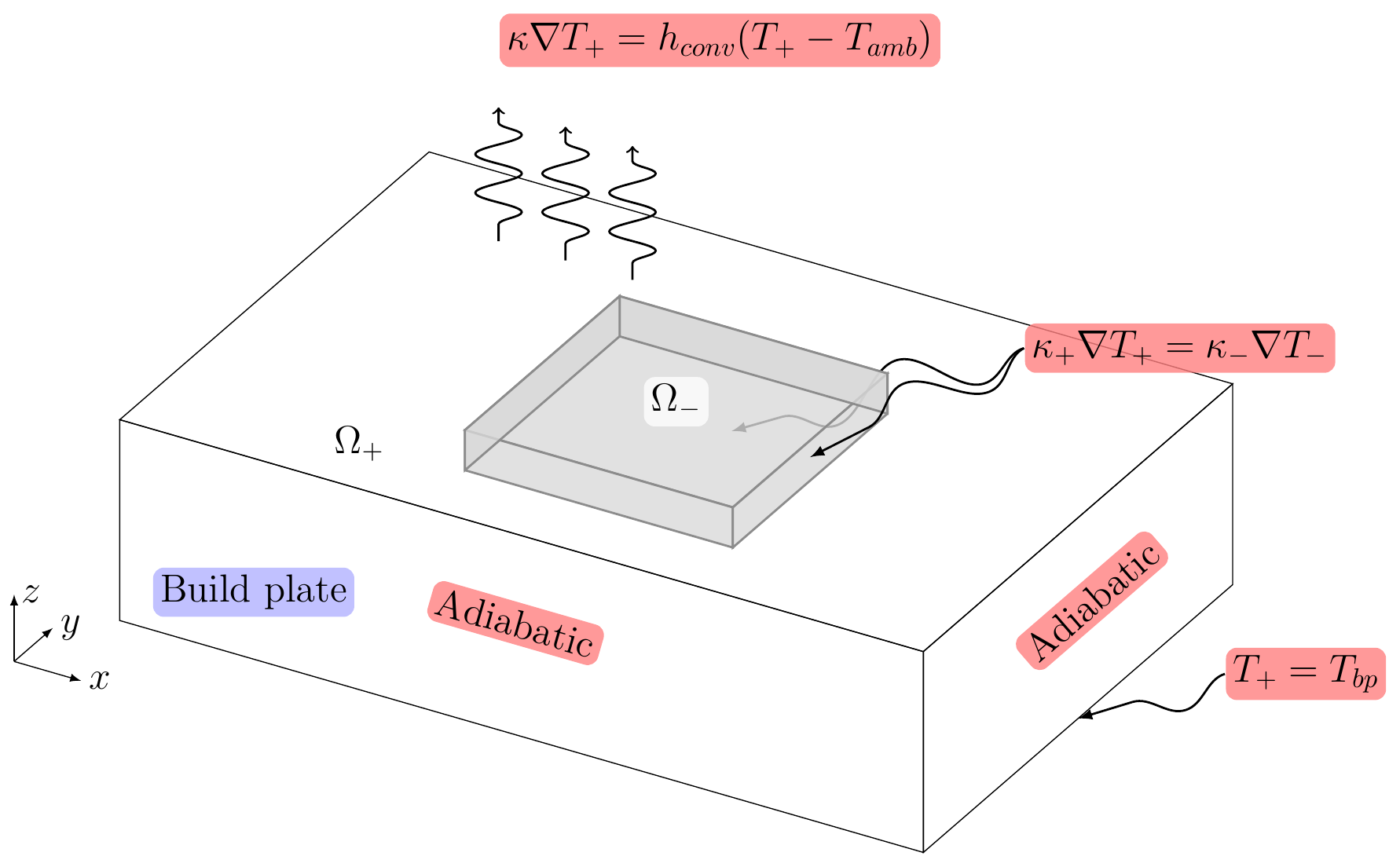}
	} 
	\\
	\subfloat[Local domain BCs.\label{subfig:localBCs}]
	{
		\includegraphics[width=0.75\textwidth]{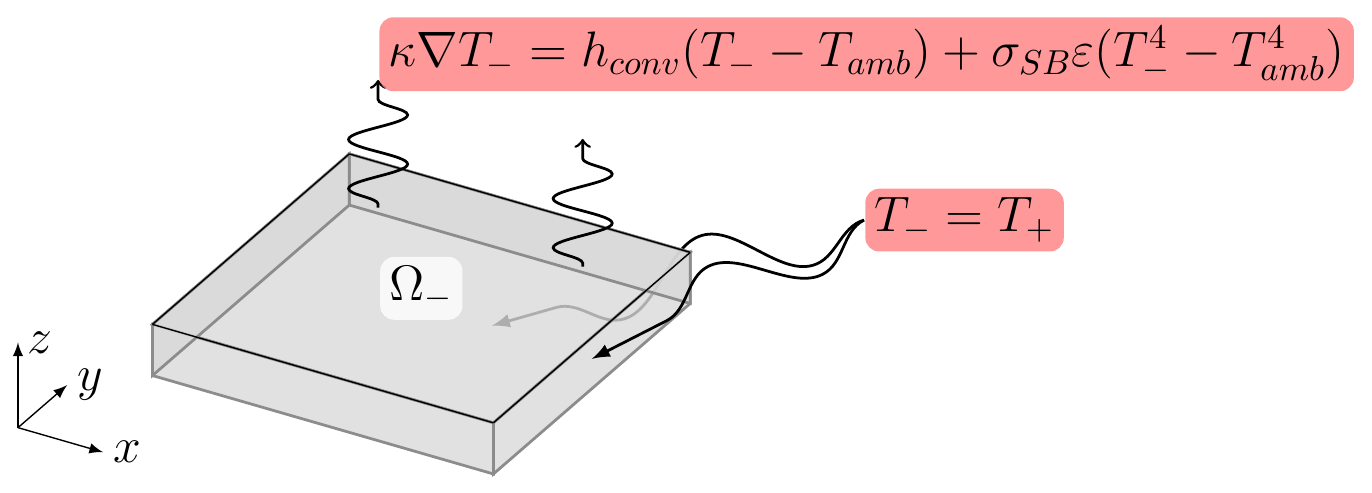}
	}
  
	\caption{Thermal problem boundary conditions (BCs). On the lower surface of the global domain (a) the temperature is fixed to the build plate temperature $T_{bp}$; adiabatic BCs are applied on the lateral surfaces of the build plate, while convection fluxes are imposed on the upper surface; finally, a zero-jump condition is imposed on the fluxes on the immersed local domain boundaries.
	In the local domain (b) the temperature continuity with the global solution is enforced by means of Dirichlet BCs on the lateral and lower surfaces, whereas heat is dissipated through a convection and radiation heat flux on the upper surface.   \label{fig:thermalBCs}}
\end{figure}
To capture the multi-scale nature of LPBF processes, the two-level method distinguishes among a local $\Omega_-$ and a global $\Omega_+$ domain, which we model differently in the proposed thermal model. \cref{fig:thermalBCs} depicts the boundary conditions (BCs) on the local and global domain. As can be observed, the local domain dissipates heat by convection and radiation through the upper surface (see \cref{subfig:localBCs}) by means of a heat loss flux term defined as follows:
\begin{equation}
	\kappa\nabla T_-= h_{conv}(T_- - T_{amb})+\sigma_{SB}\varepsilon (T_-^4 - T_{amb}^4),
\end{equation}
where $h_{conv}$ is the heat transfer coefficient by convection due to the inert gas (Argon) flow present in the chamber, $T_{amb}$ is the ambient temperature of the building chamber, $\sigma_{SB}=5.87\times 10^{-8}$ [W/m$^2$/K$^4$] is the Stefan-Boltzmann constant, and $\varepsilon$ the emissivity of the powder bed.
In the present contribution, we assume that the global domain does not dissipate heat by conduction through the lateral surfaces of the powder bed (adiabatic condition), whereas it dissipates heat by convection through the upper surface. The latter is modeled by means of a  heat loss flux term defined as:
\begin{align}
	\kappa\nabla T_+ &= h_{conv}(T_+ - T_{amb}).
\end{align}

Quantities such as emissivity and heat transfer coefficient are unknown and not readily measurable, thus they required a thorough calibration procedure (see, e.g., \cite{kollmannsberger2019}).
A further simplification we adopt in our model is that radiation effects are neglected in the global domain, assuming that they play a minor role in regions far from the heat affected region, which we consider lying entirely within the local domain. As described in \cref{subfig:globalBCs}, global domain BCs are completed by Dirichlet BCs on the lower surface of the build plate. Finally, to avoid jumps among the local and the global domain temperature fields, a zero-jump temperature condition is imposed in the local problem on the lateral and bottom surfaces of $\Omega_-$, while, in the global problem, a zero-jump temperature flux condition is enforced on the lateral and bottom surfaces of $\Omega_-$ which are immersed in the global domain.

Due to the specific features of the two-level method as presented in \cite{viguerie2020fat}, the powder-bed can be directly included in our thermal model to impose temperature and flux continuity on the immersed local domain boundaries.

\section{The two-level method in time}
While the two-level method as defined so far provides an effective way of decoupling the different spatial scales of the LPBF problem, it does not address multiscale behavior in time. Indeed, in \cite{viguerie2021numerical}, the authors used identical temporal discretizations for both local and global problems. In reality, such a uniform temporal discretization is inefficient, as a temporal discretization fine enough to resolve the local dynamics is often not necessary to give a satisfactory solution in the global domain.  
\par
\begin{figure*}
\centering
\includegraphics[width=\textwidth]{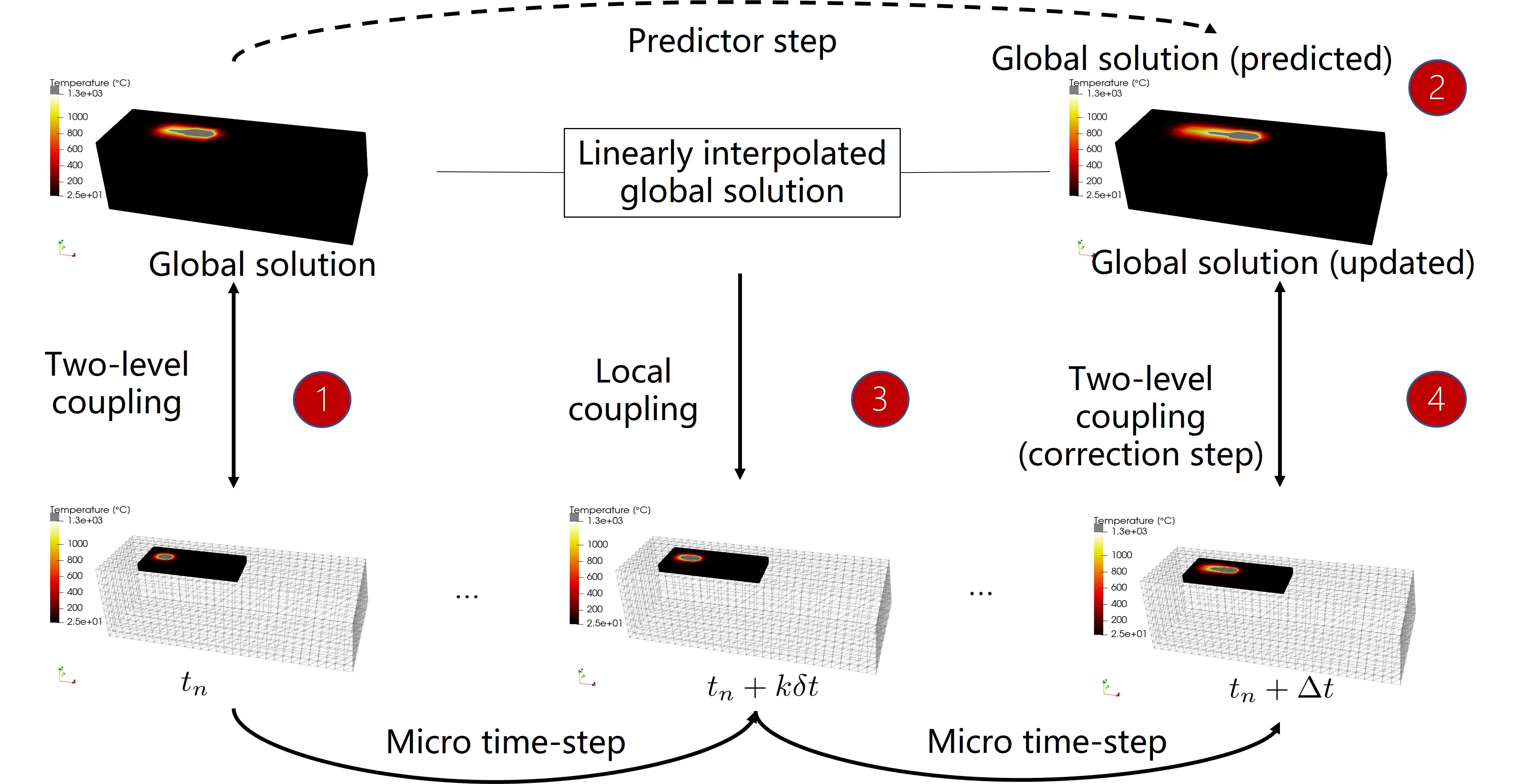} \caption{Flow diagram of the spatiotemporal two-level method: at each macro step, we solve both global and local problems using the two-level coupling approach (1). Additionally, we compute a \textit{predictor} solution of the global problem at next time step (2) such that, at each micro step, we solve a local problem using as local domain boundary conditions the global solution values obtained by means of a linear interpolation among the global solution at time $t_n$ and the predictor solution at time $t_n + \Delta t$ (3). Finally, we compute the updated global solution with the two-level coupling at time $t_n + \Delta t$ (4). }\label{fig:TwoLevelSTFig}
\end{figure*}

 We thus propose to extend the two-level in space framework into a spatiotemporal two-level method following the approach written in Algorithm 1 and detailed in Fig. \ref{fig:TwoLevelSTFig}.
\vspace{5mm}

\begin{algorithm}[h]\label{alg:TwoLevelST}
\SetAlgoLined
\KwResult{Spatiotemporal two-level method}
 At $t=0$: Solve initial two-level problem for local and global domains; \\
 \For{$n$=1:N-1}{
 At $t=t_n + \Delta t $, compute an intermediate global temperature $\widetilde{T}^+$ (predictor solution).\\
 \For{$i$=0:m-1}{
  at $t=t_n + i \delta t $ solve the time-local problem \eqref{localProblemTime};
 }
 At $t=t_n + \Delta t$: solve two-level problem for local and global domains (correction step); \\
 }
 \caption{Spatiotemporal two-level method: we solve at each one of the N macro-steps the full two-level coupled problem, as well as a predicted intermediate global temperature. At the intermediate micro-steps between each macro-step, we solve only the local problem, using a linear combination of the recently computed global temperature and predicted global temperature.}
\end{algorithm}
\vspace{5mm}
We refer to the time steps where the global solution is solved as macro steps, with time step size $\Delta t$, each separated by $m$ micro steps, with time step size $\delta t$. At each macro step, we solve the full two-level method in space at a time $t$ for a global temperature $T^+$. We then solve a \textit{predictor} step for an intermediate global temperature $\widetilde{T}^+$ at time $t+\Delta t$. 
\par For $i=0,\,1,\,2,\,...,\,m-1$, with $m$ equal to the number of micro-steps per macro-step, we then solve the local problem:

 \begin{equation}\label{localProblemTime} \begin{alignedat}{2}
 	&c_-\rho_- \frac{\partial T^-}{\partial t} + \chi \rho_- \frac{\partial f_{pc}}{\partial t} - \nabla \cdot \left(\kappa_- \nabla T^- \right) &&\\ &= Q \qquad &&\text{  in  } \Omega_- \times \lbrack t_0,\,t_{end}\rbrack \\
 	&T^- = \left(1-\frac{i}{m}\right)T^+ + \frac{i}{m} \widetilde{T}^+ \qquad &&\text{  on  } \gamma_D \\
 	&\kappa_- \frac{\partial T^- }{\partial\boldsymbol{n} } = 0 \qquad &&\text{  on  } \Gamma_N .
\end{alignedat}\end{equation}

 Finally, at the next macro step, we use the recently computed local solution to obtain the global solution $T^+$ at $t+\Delta t$ using the standard two-level iteration scheme described in \cref{ssec:TwoLevelMethod}, obtaining then a corrected version of $\widetilde{T}^+$. In this way, the algorithm can be viewed as a predictor-corrector type method. We note that a specific choice of time integration scheme is not specified. In principle, the presented framework is independent of such a selection; the two-level time integration scheme may be applied in combination with different numerical integration methods.
\par Combining the two-level method in space with a multirate time integration scheme as shown above offers several notable advantages. The geometric flexibility that makes the method attractive for spatial integration also has similar benefits in the temporal domain. The degrees of freedom can be moved and changed freely without requiring changes to the mesh topologies (as is the case with refinement and derefinement). When compared to a method as the one shown in \cite{hodge2021towards}, there is no need to track the `fast' and `slow' degrees of freedom specifically, as one can use the different meshes to directly accomplish this task. Lastly, the fact that this method allows for the easy application of uniform meshes allows for many advantages in the temporal domain, as CFL-type conditions and other stability considerations become easier to interpret.
\par There are many theoretical questions with the proposed algorithmic framework. Namely, we would like to know whether convergence orders are maintained, stability conditions, and how refinements to the micro- and macro-time scales affect error behavior. Similar investigations were carried out for the two-level method in space \cite{viguerie2020fat}. In the following section, we will address some of these concerns, notably the dependence of error behavior on micro- and macro-time-step size empirically via numerical experiment. However, we intend to perform a more comprehensive investigation and formal analysis of such accuracy, as well as stability, characteristics in future work. 
\section{Numerical Experiments}\label{sec:experimentalValidation}
All of the two-dimensional results presented in this section have been obtained on a 2020 Apple MacBook Pro with 2 GHz quad-core Intel i5 processor and 16 Gb RAM. For the three dimensional results, we deployed an Ubuntu 20.4 OS equipped with an Intel\textsuperscript{\tiny\textregistered} Xeon\textsuperscript{\tiny\textregistered} W-2125 CPU @ 4.0GHz and 256 Gb RAM. All the numerical schemes, in both two- and three-dimensions, have been implemented in the FreeFem++ environment, a partial differential equation solver written in C++ \cite{MR3043640}.

\subsection{Two-Dimensional Convergence Study}
\begin{figure}
\centering
\includegraphics[width=0.5\textwidth]{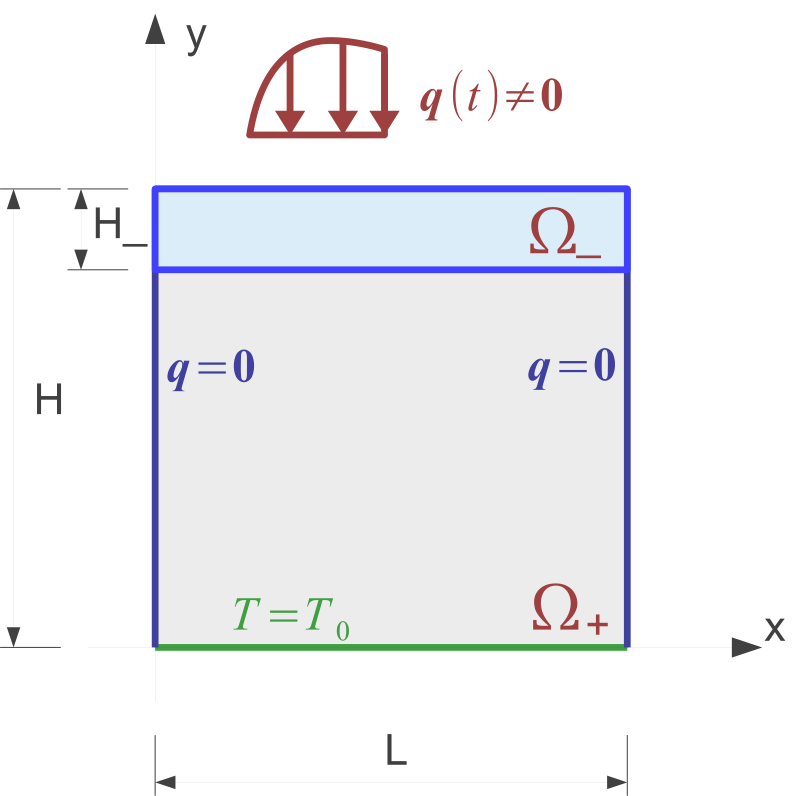}\caption{Problem setup for two-dimensional convergence study.}\label{fig:TwoLevel2DSetup}
\end{figure}

We first examine the error behavior on a two dimensional model problem. We consider the thermal model \eqref{thermalEqn}, with the material IN625 (thermal parameters defined in \cite{Mills2002}). In order to evaluate the convergence, we compare the two-level method to a reference solution resolved monolithically in space and time. The latent heat term is considered only in the local domain. As the spatial error of the two-level method was already considered in \cite{viguerie2020fat, viguerie2021numerical}, we consider identical spatial discretizations for the reference, local, and global problems in order to analyze the temporal error exclusively.
\par We consider a general problem setup similar to the problems shown in \cite{viguerie2020fat, viguerie2021numerical} and depicted in Fig. \ref{fig:TwoLevel2DSetup}. We set $L$= 5 mm, $H$ = 1 mm, $H_-$ = 0.375 mm. The laser term $Q$ is defined as:
\begin{align}\label{QDefn2D}
Q(t) &= P \frac{\eta}{r d} \exp \left(-\frac{(x-x_c)^2}{r^2}-\frac{(y-y_c)^2}{d^2}\right).
\end{align}
We set the laser power $P$=1.8 W, the absoptivity of the material $\eta=1$, laser depth $d$ = 0.0125 mm, laser radius $r = 0.1$ mm, and laser height $y_c=1$ mm. The laser center $x_c(0)=0.5$ mm and moves with the speed 0.01 mm/s. We define the phase change function as:
\begin{align}\label{2dPhaseChange}
f_{pc} &= \frac{1}{2} S \left\lbrack 1 - \tanh \left(S-\frac{(T_s+T_l)}{2} \right)^2\right\rbrack
\end{align}
where $T_s$ is the solid temperature and is set to 1290$^\circ$ C, while $T_l$ is 1380$^\circ$ C \cite{Mills2002}. The parameter $S$, a dimensionless quantity corresponding to the sigmoid sharpness, is set to $S=0.05$.
\par We consider two one-second sweeps of the laser across the geometry, and compare the error in the local domain at each time-step to the reference monolithic solution. 
\par We first fix a local (and reference) time-step size of 0.01 s and varying global time step sizes of 0.2, 0.1, 0.05, and 0.02 s. We seek to analyze the error behavior as the global time-step size is reduced in relation to the local time-step size. We note that the difference with the reference solution is exact up to convergence tolerance for global and local time-step size identical to the reference time-step size.  
\begin{figure}
\centering
\includegraphics[width=\textwidth]{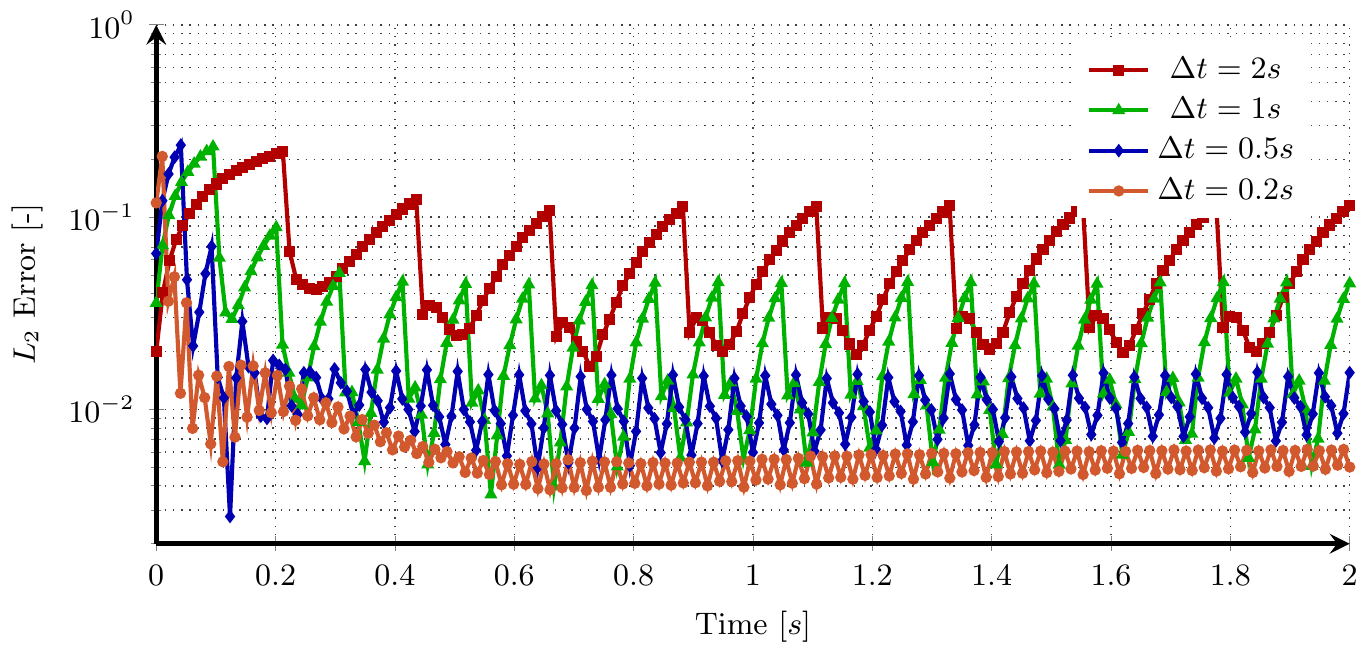}\caption{Error behavior in time, two-dimensional convergence study. We see that, for fixed local time-step size, reducing the global time-step size reduces the error.}\label{fig:TwoLevelTemporalErrorFig}
\end{figure}

\par In Fig. \ref{fig:TwoLevelTemporalErrorFig}, we plot the error at each time-step compared to the reference solution for each global time-step size. We observe a general reduction in both overall error magnitude and error fluctuation as we reduce the global time-step size. An important observation is the qualitative error behavior as time evolves. We notice a fluctuating error trend; this is caused by the local-global coupling in time. In the first local time-steps, immediately after a global time-step, the error is lower. This error then increases progressively as the distance in time from the last global solution increases. Upon the solution of a new global problem, this error then drops once again, following the expected behavior, as the global solution used to obtain local boundary conditions becomes a worse approximation to the true boundary conditions as time evolves. Naturally, the magnitude of these variations is related to the chosen global time-step size. 
\par
\begin{figure}
\centering
\includegraphics[width=\textwidth]{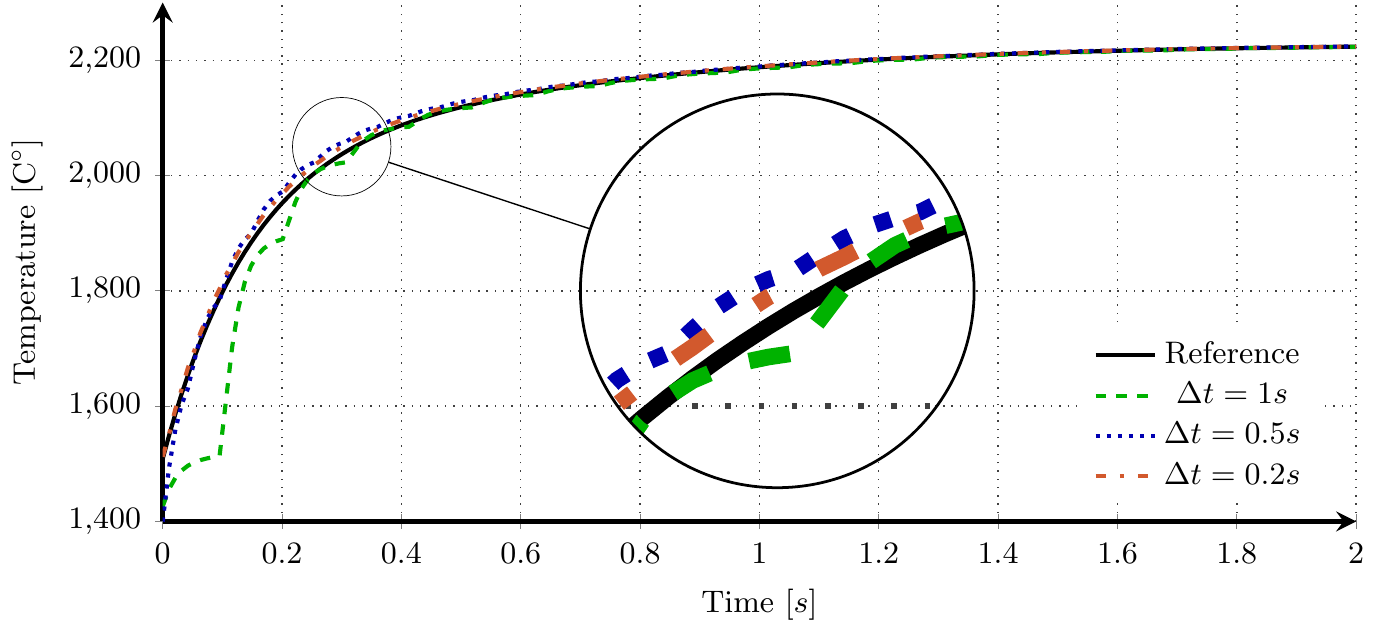}\caption{Temperature integrated in space over the line connecting (0,0.99),(5,0.99). We see a closer match to the physical temperature as $\Delta t_{Glob}$ decreases. For all global time-step sizes, the two-level temperature better matches the reference temperature as the time increases. }\label{fig:TwoLevelTempMatch}
\end{figure}

 Fig. \ref{fig:TwoLevelTempMatch} reinforces these findings, and shows the evolution of the quantity:
\begin{align}\label{controlTemp}
T^* &= \int_0^L T(x,\,0.99)\,dx    
\end{align}
in time. As expected, we see that, as the global time step decreases, the curves match the reference solution more and more closely. In general, we also see that, for all global time-step levels, the majority of the error occurs at the beginning of the simulation, and the solutions improve in quality as time evolves. For the case in which $\Delta t_{Glob} = 0.1$s, we very clearly observe the effect of the solution being accurate near the solutions of each global problem, with the discrepancy between the two-level and reference solutions increasing as the distance from the global solution increases.
\par Having analyzed the error behavior for a fixed local time-step and varying global time-step, we then seek to observe the opposite relationship. We run simulations for global time steps of 0.2 s, 0.1 s, and 0.05 s. Then, for each global time step, we consider local time steps of 0.050 s, .025 s, and 0.010 s (note for the case of $\Delta t = .05$, we consider only $\delta t=0.025,\,.010$ s). We then quantify the error behavior in terms of mean relative $L^2$ error over the entire time interval as compared to a reference solution (computed with a time step of .0025 seconds). 
\par
\begin{figure}
\centering
\includegraphics[width=\textwidth]{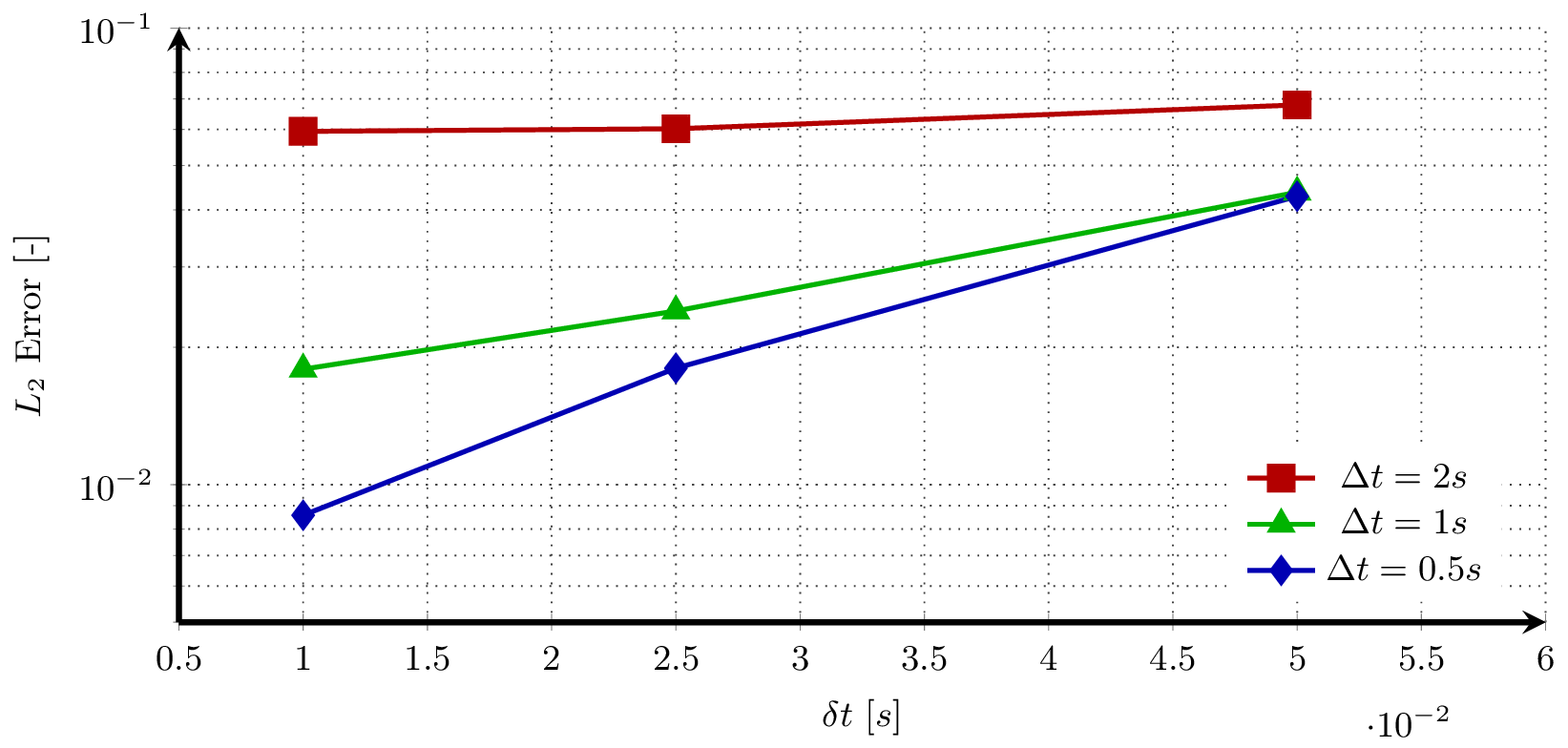}\caption{Mean error compared to a reference solution for a range of local and global solutions. For each global time-step level, we observe an improved error behavior from refining the local mesh. The global error also plays a significant role, and we see that, after a certain point, refinements of the local time-step without further global time-step refinements yield diminishing returns.  }\label{fig:ErrorPlotGlobalLocal}
\end{figure}

 We plot the relevant results in Fig. \ref{fig:ErrorPlotGlobalLocal}. and again observe the general expected behavior; a lower global time-step size reduces error, with refinements in the local time-step size also further reducing error. We see error reduction independently of the global time-step size; however, the degree in which the micro time step size decreases the error appears closely related to the macro time step size. For the large $\Delta t$ of 0.2 s, we see a very slight effect of time-step refinement on error behavior. This effect becomes more pronounced for $\Delta t$=0.1 s, and even more so for $\Delta t$=0.05 s. This suggests that, in order for $\delta t$ to have a significant effect on the error behavior, we require a sufficiently small $\Delta t$. This is indeed similar to the behavior observed for the spatial error in \cite{viguerie2020fat, viguerie2021numerical}, in which it was shown that both mesh resolutions influence error in an interdependent way, with finer global mesh resolutions leading to a higher influence of local-scale error. A more rigorous formal analysis of this error behavior and the dependence $\Delta t$ and $\delta t$ will be the subject of a future work.

\subsection{Three-dimensional single track example}
The second numerical test is a simulation of a single-track laser in three dimensions. This test features a fixed global mesh and a local mesh that moves in time, following the laser. We seek to observe the following important aspects of the method:
\begin{enumerate}
    \item[i)] The advantage granted in terms of computational time by applying the two-level method in time over a space-only two-level method with fine temporal refinement across all scales;
        \item[ii)] The possibility of moving the local mesh in time, providing an efficient and simple way to separate the different temporal scales in the problem, without enduring a loss in accuracy.
\end{enumerate}
\par Also in this 3D numerical example, we consider a laser heat source moving along a straight line of 2 mm along the $x-$axis on the upper surface of a bare plate of IN625. In this case, we set $L$= 3mm. $H$=1mm, $H_-$ =0.1 mm. For a laser located at a given instant in time in ($x_c, y_c, z_c$), we model the laser heat source term using a Gaussian heat source $Q_-$ in the local domain and a distributed constant heat source model $Q_+$ in $\Omega_h\subset\Omega_+$ defined respectively as:
\begin{align}\begin{split}
Q_-(t) &=  \frac{6 \sqrt{3} P \eta}{2\pi r^2 d}  \\
&\exp \left(-\frac{3(x-x_c)^2}{r^2} - \frac{3(y-y_c)^2}{r^2} -\frac{3(z-z_c)^2}{d^2} \right)\end{split}\label{QDefn2D}\\
Q_+(t) &=  \frac{P \eta}{r d v \Delta t}
\end{align}
with $\Omega_h$ having dimensions $v\Delta t \times r\times d$.
Such a definition of $Q_+$ allows us to adopt larger time steps in the global problem without obtaining a jumping behavior of the global solution, while at the same time the accuracy in proximity of the laser source is kept by the Gaussian heat source $Q_-$ in the local domain.
\par We set the laser power $P$=179.2 W, the laser depth $d$=0.050 mm, the laser radius $r=0.085$ mm, the absorptivity $\eta = 0.38$, and the laser speed $v=800$ mm/s. We define the phase change function $f_{pc}$ as in the previous 2D example.
\par To evaluate the accuracy and the efficiency of the proposed spatiotemporal two-level approach, we first compute a reference solution of the problem employing a uniform linear tetrahedral mesh with element length size $h=0.02$ mm and a uniform time step $\Delta t = 1.e-4$ s. We also consider a spatial two-level approach with a moving local domain of $1.2\times 0.5 \times 0.1$ mm$^3$, a local element length size $h_-=0.02$ mm, a global element length size $h_+=0.1$ mm, and a uniform time step $\Delta t = 10^{-4}$ s. Finally, we adopt the proposed spatiotemporal two-level scheme using the same spatial discretization of the spatial two-level  model but with a local time step $\delta t=10^{-4}$ and a four time coarser global time step $\Delta t=4\cdot 10^{-4}$.
\par
\begin{figure}[h!]
	\centering
	\subfloat[Local temperature distribution at t=0.0011.\label{subfig:moving11}]
	{
		\includegraphics[width=0.75\textwidth]{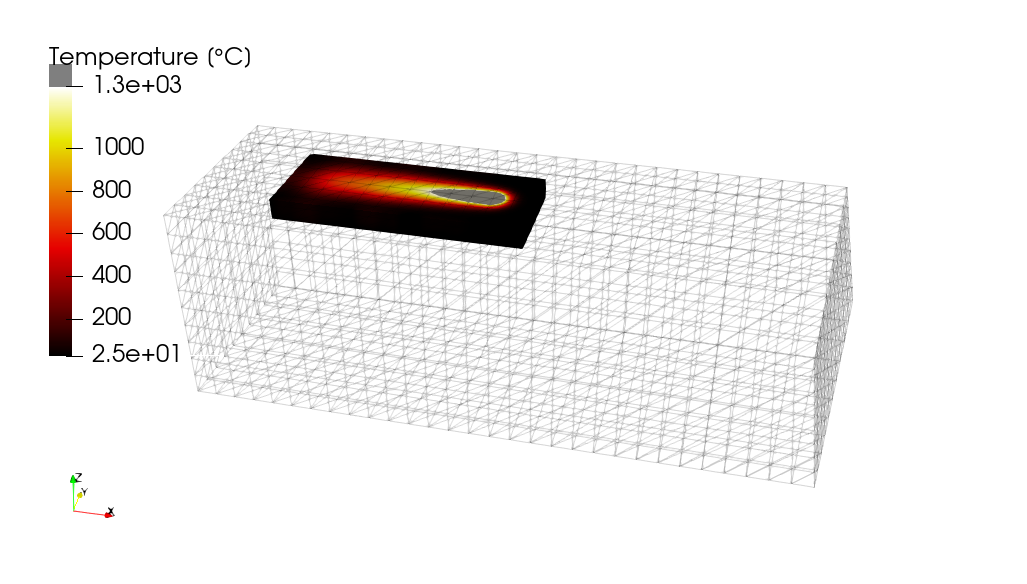}
	}
	\\
	\subfloat[Local temperature distribution at t=0.0012.\label{subfig:moving12}]
	{
		\includegraphics[width=0.75\textwidth]{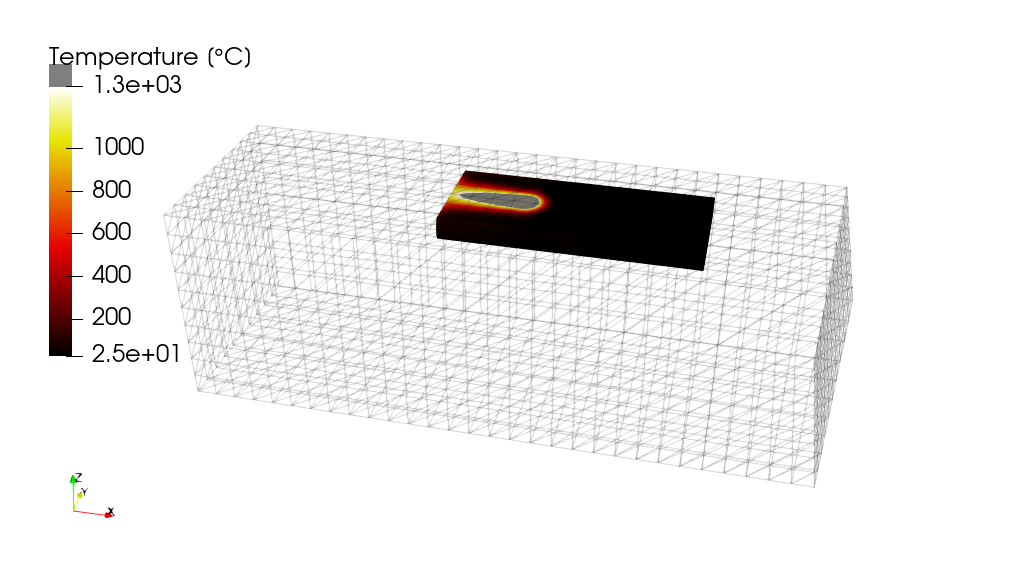}
	}    
  \caption{Moving local mesh temperature distributions.\label{fig:movingLocal}}
\end{figure}

 \Cref{fig:movingLocal} shows the temperature distribution in the local domain at two consecutive time steps before and after a movement of the local domain. Thanks to the adopted two-level coupling scheme described in \cref{ssec:TwoLevelMethod}, no projection is required during this step since the global solution carries the steady-state solution which is transferred to the new local domain through the domain boundary coupling.

%
%
\par
\begin{figure*}
\centering
\includegraphics[width=\textwidth]{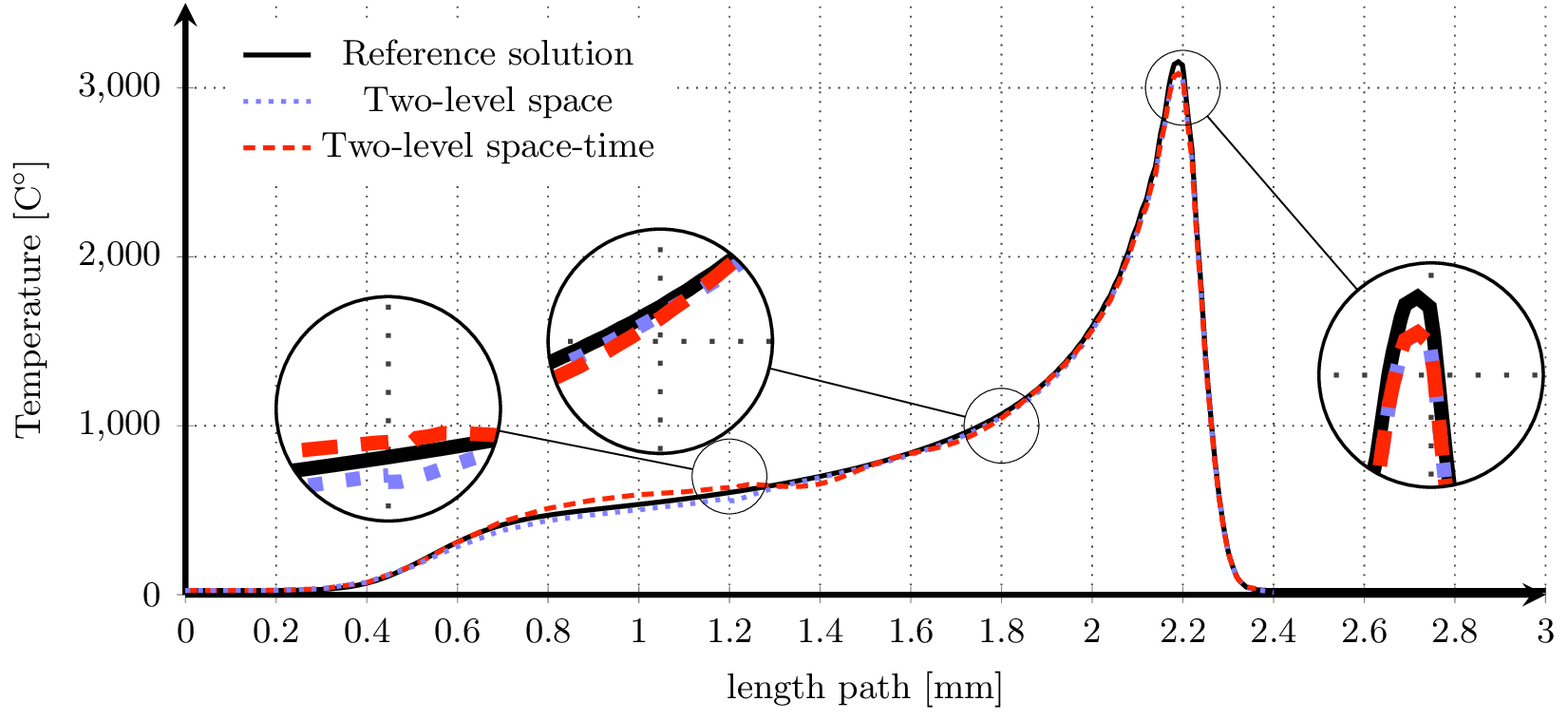} \caption{Temperature distribution along the laser moving direction at time $t=0.02$ s. The reference solution is obtained using a uniform mesh and monolithic time integration, the two-level space solution using a moving local mesh and monolithic time integration, finally the two-level space-time employs the presented two-level spatiotemporal approach. No significant differences can be observed among the three solutions. }\label{fig:TwoLevelTempPlot}
\end{figure*}

 In \cref{fig:TwoLevelTempPlot}, the temperature distribution at time $t=0.020$ s, i.e., after 20 micro time-step, is depicted for all the three different spatiotemporal discretizations defined above. As can be noticed, both the spatial two-level and the spatiotemporal two-level apporaches are able to closely capture the reference temperature solution, making us confident on the accuracy of the proposed numerical framework.
\subsection{Multi-track alternate scan path example}
In this last numerical example, a spatiotemporal two-level simulation is computed for a laser source traveling on a $6\times 6\times 2$ mm$^3$ global domain following 50 scan tracks with alternate scan directions. The local domain size is set to $2\times 0.5\times 0.2$ mm$^3$, while the process parameters the same as in the previous example.
\begin{figure*}[h!]
	\centering
	\subfloat[t=0.0045 s.\label{subfig:startpoint}]
	{
		\includegraphics[width=0.5\textwidth]{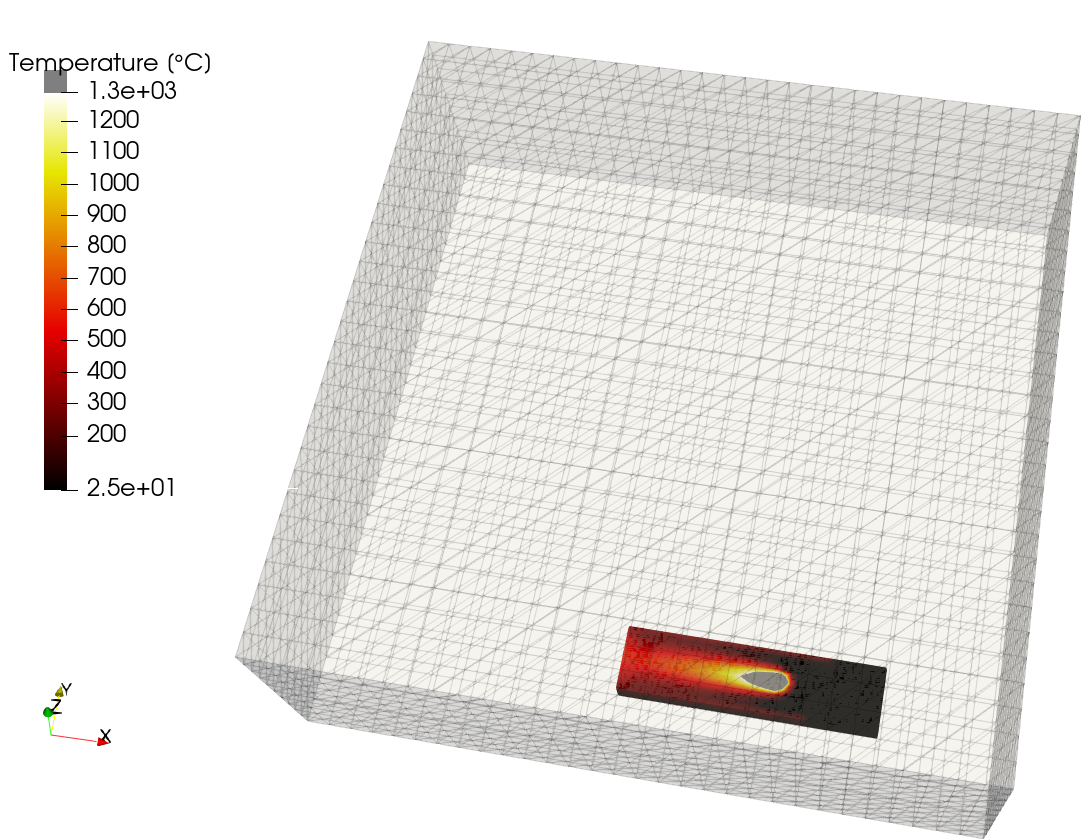}
	}
\hspace{50pt}
	\subfloat[t=0.1 s.\label{subfig:turningpoint}]
	{
		\includegraphics[width=0.4\textwidth]{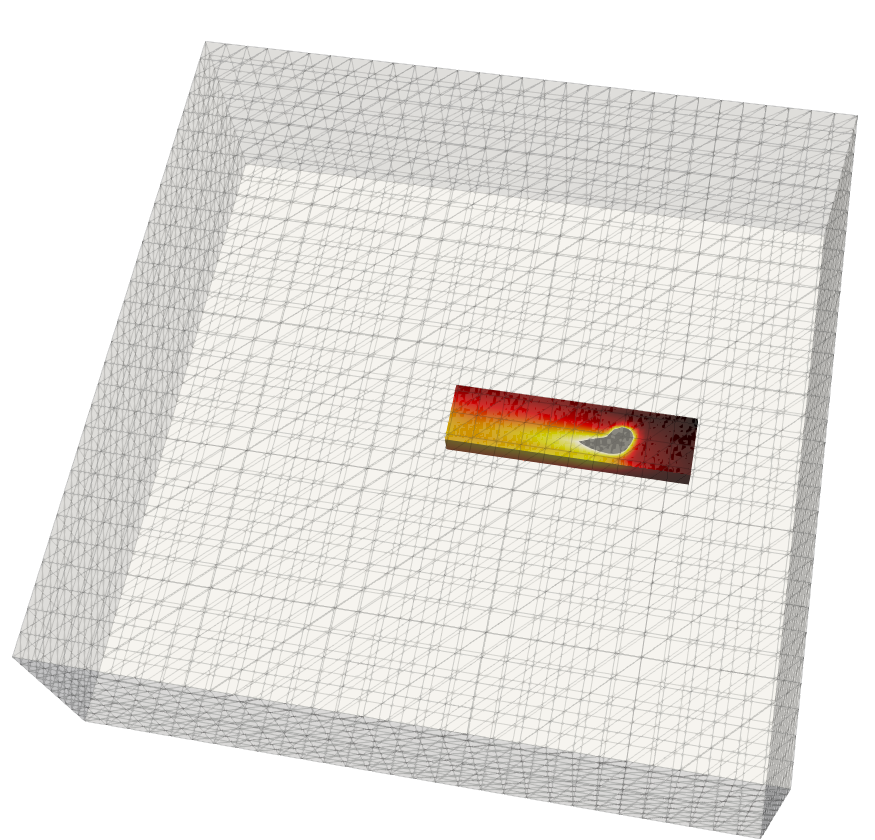}
	}  
\hspace{50pt}  
		\subfloat[t=0.21595 s.\label{subfig:endpoint}]
	{
		\includegraphics[width=0.425\textwidth]{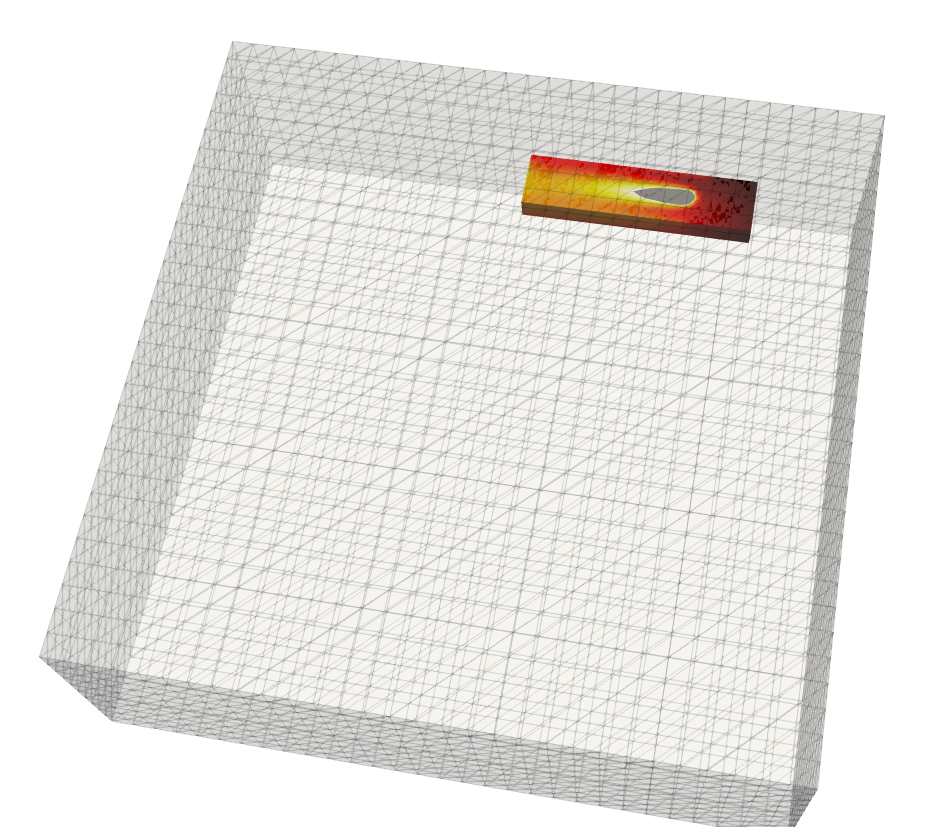}
	} 
  \caption{Melt pool shape and temperature distribution in the local domain at three different times.\label{fig:movingMulti}}
\end{figure*}
\Cref{fig:movingMulti} shows the temperature distribution and the melt pool shape at three different instants in time. Comparing \cref{subfig:startpoint} and \cref{subfig:endpoint}, it is possible to observe the influence of the residual heat of the previously scanned regions on the temperature distribution and the melt pool morphology. Such an effect is captured thanks to the two-level coupling between the local and the global domain. Moreover, also at the turning point of the scan path (see \cref{subfig:turningpoint}), the two-level coupling efficiently captures the complex temperature field distribution. 
\par
\begin{figure}[h!]
\includegraphics[width=.8\textwidth]{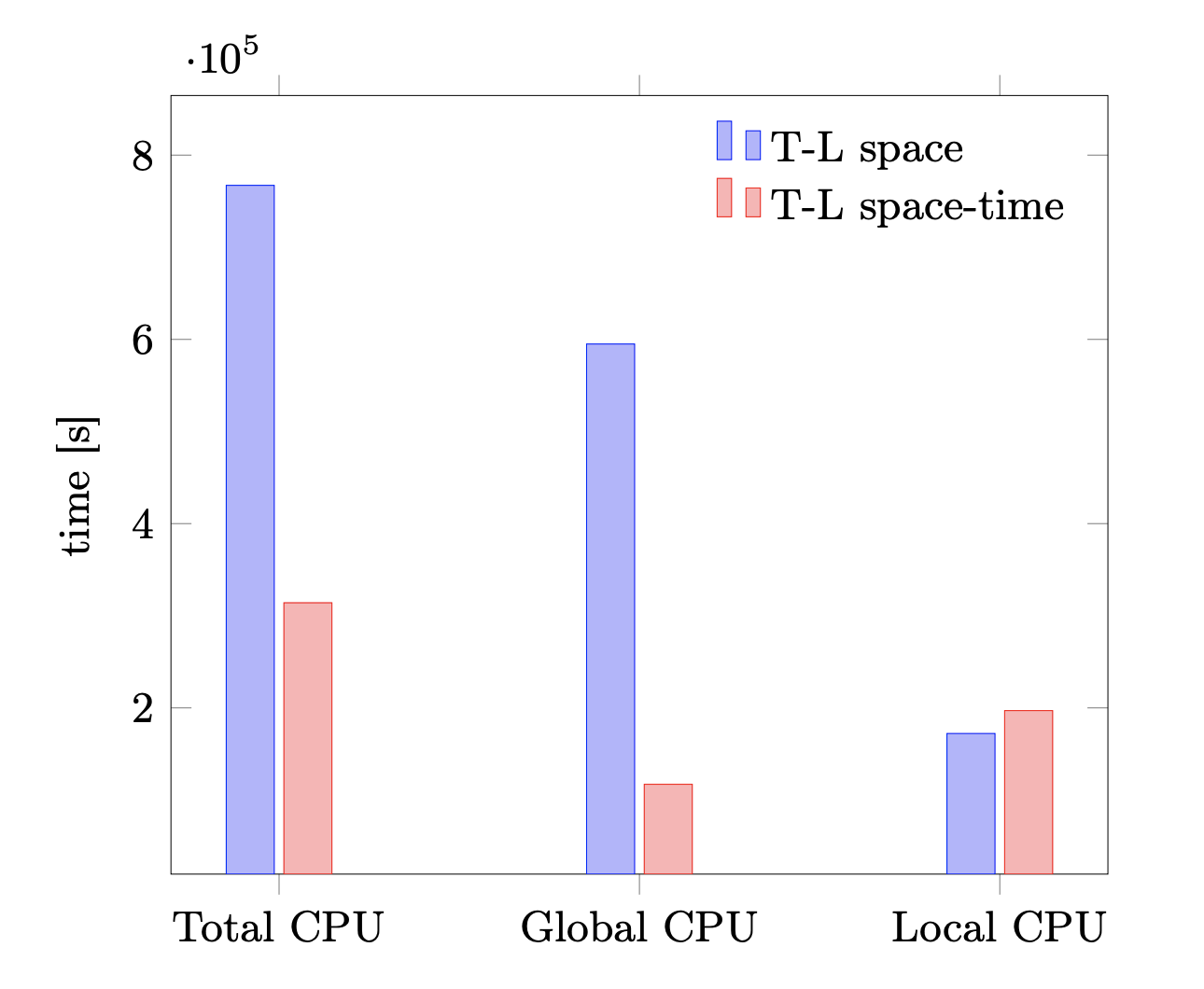}
    \caption{CPU time for different discretization schemes}
    \label{fig:TimeComparison}
\end{figure}

 Finally, such an analysis has been computed also employing the purely spatial two-level approach, setting both the global and the local time-step size equal to the local time-step size of the spatiotemporal analysis. 
\Cref{fig:TimeComparison} depicts the total CPU time, the CPU time spent to solve the global problems, and the CPU time spent to solve the local problems using both the spatial and the spatiotemporal two-level method. On the one hand, the time spent to solve the local problem remains almost constant and the difference among the spatial and spatiotemporal version of the two-level approach is primarily due to faster convergence of the former. On the other hand, we obtain a speed-up of approximately $\times 2.44$ in the solution of the overall problem. Such a computational speed-up is achieved thanks to the adopted multi-rate time integration scheme which allows us to employ a much coarser time step size for the global problem without loosing much accuracy in the solution. Moreover, the proposed numerical framework - due to its structured, non-conform nature - is embarrassingly parallelizable, thus further code optimization will allow us to exploit the effects of the adopted spatiotemporal approach even further.

	
		]%


\section{Conclusions}\label{sec:conclusions}
We extended the two-level framework first introduced in \cite{viguerie2020fat, viguerie2021numerical} for multiple spatial scales by introducing a multi-rate time integration scheme capable of resolving multiple scales in time. Such an extension is important for LPBF AM problems, as a large separation of scales occurs at both the spatial and temporal levels. We proceeded to perform an empirical error analysis, demonstrating the influence of both macro and micro time scale resolution on overall error behavior. We then performed a series of three-dimensional experiments on a realistic problem configuration. These experiments in 3D established that one may employ the spatiotemporal two-level method without incurring significant losses in accuracy, and, further, the two-level method can be combined with a moving fine-scale mesh while retaining a high-quality solution. This endows the proposed method with remarkable geometric flexibility, a crucial feature for effective LPBF solvers. It was also shown that employing the spatiotemporal two-level approach, instead of the purely spatial one, offers the potential for significant savings in computational cost by reducing the number of necessary global problem solutions.
\par There are several important directions for the present work. As this work was primarily a proof-of-concept, the error analysis was restricted to an empirical study. However, a more rigorous formal analysis of the temporal error, and how it depends on the macro and micro time-step sizes, is important for a more comprehensive understanding of the method performance. Extending the current framework from the temperature-only problem shown here, into a multiphysics setting, considering also mechanical behavior, fluid dynamics, and 3-phase compartment modelling, is also an important extension. Finally, from the computational performance point-of-view, code optimization featuring improved preconditioning and parallelization schemes is also important to fully maximize the potential computational advantages of the discussed method.

\section*{Acknowledgments}
This work was partially supported by the Italian Minister of University and Research through the MIUR-PRIN projects "A BRIDGE TO THE FUTURE: Computational methods, innovative applications, experimental validations of new materials and technologies” (No. 2017L7X3CS) and "XFAST-SIMS" (no. 20173C478N).
\bibliographystyle{elsarticle-num-names}
\biboptions{sort&compress}
\bibliography{TwoLevel}

\end{document}